\documentclass[12pt]{amsart}

\usepackage{a4wide}
\usepackage{pxfonts}
\usepackage{wasysym}
\usepackage{amssymb,amsxtra,mathrsfs}
\usepackage{eucal}
\usepackage{setspace,upgreek}
\usepackage{mathtools}
\usepackage{url}
\newcommand{\mf}[1]{\mathfrak{#1}}
\newcommand{\mc}[1]{\mathcal{#1}}

\newcommand{\mr}[1]{\mathrm{#1}}
\newcommand{\lr}[2]{\langle #1,#2\rangle}
\newcommand{\up}{\upharpoonright}

\newcommand{\un}{\mathbf{1}} 
\newcommand{\N}{\mathbb{N}}
\newcommand{\R}{\mathbb{R}}
\newcommand{\C}{\mathbb{C}}
\newtheorem{theorem}{Theorem}
\newtheorem{definition}[theorem]{Definition}
\newtheorem{lemma}[theorem]{Lemma}

\newtheorem{corollary}[theorem]{Corollary}

\theoremstyle{definition}

\bibliography{database}
\begin{document}
\title[Construction of strong derivable maps]
{Construction of strong derivable maps via functional calculus of unbounded spectral operators in Banach spaces}
\author[B. Silvestri]{Benedetto Silvestri}

\address{Benedetto Silvestri\\
440 W Kelly Ave Apt 6\\
Jackson, WY 83001-8544\\
United States\\
(via express couriers)}

\address{
Benedetto Silvestri\\
PO BOX 14845\\
Jackson, WY 83002-4845\\
United States\\
(via United States Postal Service)}

\email{BNDSILVESTRI@PROTON.ME}

\date{\today}
\keywords{strong operator derivability, derivability of locally convex space valued maps, unbounded spectral operators in Banach
spaces, functional calculus}
\subjclass[2020]{46G05, 47B40, 47A60}
\begin{abstract}
We provide sufficient conditions for the existence of a strong derivable map and calculate its derivative 
by employing a result in our previous work on strong derivability of maps arising by functional calculus of 
an unbounded scalar type spectral operator $R$ in a Banach space and the generalization to complete locally convex spaces of a 
classical result valid in the Banach space context. We apply this result to obtain a sequence of integrals 
converging to an integral of a complete locally convex space extension of a map arising by functional calculus of $R$. 
\end{abstract}
\maketitle
%%%%%%%%%%%%%%%%%%%%%%%%%%%%%%%%%%%%%%%%%%%%%%%%%%%%%
%%%%%%%%%%%%%%%%%%%%%%%%%%%%%%%%%%%%%%%%%%%%%%%%%%%%%
%%%%%%%%%%%%%%%%%%%%%%%%%%%%%%%%%%%%%%%%%%%%%%%%%%%%%
\section*{Introduction}
Given a sequence $\{f^{n}\}_{n\in\N}$ of suitable analytic functions defined on an open neighbourhood $U$ of the spectrum of a 
possibly unbounded scalar type operator $R$ in a complex Banach space, an application of our result \cite[Theorem 1.23]{sil0} 
establishes that for every $v$ in the domain of $R$ and $n\in\N$, the map $W\ni t\mapsto f^{n}(tR)v$ is derivable in
$W$ an open interval of $\R$ and its derivative equals
\begin{equation*}
\frac{df^{n}(tR)v}{dt}=R\frac{df^{n}}{d\lambda}(tR)v,
\end{equation*}
where $g(R)$ is the Borelian functional calculus of $R$ for any Borelian map $g$ defined on $U$.
In the present paper we construct a locally convex space $F$ and employ the above equality to find conditions ensuring the 
existence of a sequence $\{\mf{f}^{n}\}_{n\in\N}$ in $\mc{C}^{1}(W,\widehat{F})$, a sequence 
$\{\mf{g}^{n}\}_{n\in\N}$ in $\mc{C}(W,\widehat{F})$, a map $\mf{g}^{\infty}\in\mc{C}(W,\widehat{F})$ and a subsequence 
$\{\mf{g}^{n_{k}}\}_{k\in\N}$ of $\{\mf{g}^{n}\}_{n\in\N}$ such that 
\begin{equation*}
\begin{aligned}
(\forall n\in\N)(\frac{d\mf{f}^{n}}{dt}&=\mf{g}^{n});
\\
\lim_{k\in\N}\mf{g}^{n_{k}}&=\mf{g}^{\infty}\text{ in }\mc{C}_{u}(W,\widehat{F}).
\end{aligned}
\end{equation*}
Here $\widehat{F}$ is the Hausdorff completion of $F$, $\mc{C}_{u}(W,\widehat{F})$ is the uniform space of continuous maps 
defined on $W$ and valued in $\widehat{F}$ provided by the uniformity of uniform convergence, 
$\mc{C}^{1}(W,\widehat{F})$ is the linear space of derivable maps defined on $W$ and valued in $\widehat{F}$ with continuous 
derivative. The two equalities above imply 
\begin{equation*}
\lim_{k\in\N}\frac{d\mf{f}^{n_{k}}}{dt}=\mf{g}^{\infty}\text{ in }\mc{C}_{u}(W,\widehat{F}).
\end{equation*}
Now the above limit is one of the two hypotheses of the extension of the result \cite[Corollary 8.6.4]{die1} to complete locally 
convex linear space valued maps defined on open sets of $\R$. As a consequence we establish in \textbf{Theorem \ref{05221753}} 
sufficient conditions guaranteeing the existence of a map $\mf{f}^{\infty}\in\mc{C}^{1}(W,\widehat{F})$ such that 
\begin{equation*}
\begin{aligned}
\lim_{k\in\N}\mf{f}^{n_{k}}&=\mf{f}^{\infty}\text{ in }\mc{C}_{u}(W,\widehat{F}); 
\\
\frac{d\mf{f}^{\infty}}{dt}&=\mf{g}^{\infty}.
\end{aligned}
\end{equation*}
We then apply this result in \textbf{Corollary \ref{11251048}} to obtain for every $u_{1},u_{2}\in W$ the following limit of 
integrals
\begin{equation*}
\lim_{k\in\N}\int_{u_{1}}^{u_{2}}\mf{g}^{n_{k}}(t)dt
=
\int_{u_{1}}^{u_{2}}\mf{g}^{\infty}(t)dt.
\end{equation*}
%%%%%%%%%%%%%%%%%%%%%%%%%%%%%%%%%%%%%%%%%%%%%%%%%%%%%
%%%%%%%%%%%%%%%%%%%%%%%%%%%%%%%%%%%%%%%%%%%%%%%%%%%%%
%%%%%%%%%%%%%%%%%%%%%%%%%%%%%%%%%%%%%%%%%%%%%%%%%%%%%
\section*{notation}
In this paper we follow the notation and definitions provided in \cite{sil0} which therefore will not be repeated made exception
for what follows.
\par
If $X,Y$ are topological spaces let $\mc{C}(X,Y)$ be the set of continuous maps defined on $X$ and valued in $Y$, 
if $Z$ is a uniform space, then $\mc{C}_{u}(X,Z)$ and $\mc{C}_{s}(X,Z)$ denote the uniform space of continuous maps defined on 
$X$ and valued in $Z$ provided by the uniformity of uniform convergence and simple convergence respectively.
Let $G$ be a complex Banach space, $R$ be a possibly unbounded scalar type operator in $G$ and $E$ its resolution
of the identity. If $D$ is a linear subspace of $G$, then we let $L(D,G)$ be the linear space of linear morphisms from $D$ to
the linear space underlying $G$, while $B(G)$ denotes the normed space of linear continuous so, 
bounded endomorphisms of $G$ with the usual $\sup$-norm. Let $\mc{B}(\C)$ be the set of Borelian subsets of $\C$ 
and for any $U\in\mc{B}(\C)$ let $\mr{Bor}(U)$ be the complex linear space of the complex valued Borelian maps defined on $U$.
Let $U\in\mc{B}(\C)$. For any map $p:U\to\C$ we let $\widetilde{p}:\C\to\C$ be the $0$-extension of $p$ namely
$\widetilde{p}\up U=p$ and $\widetilde{p}(\lambda)=0$ if $\lambda\in\complement_{\C}U$. For each map $F:\C\to\C$ define 
\begin{equation*}
|F|_{\infty}^{E}
\doteqdot
E-\mr{ess}\sup_{\lambda\in\C}
|F(\lambda)|
\doteqdot
\inf_{\{\delta\in\mc{B}(\C)\mid E(\delta)=\un\}}
\sup_{\lambda\in\delta}
|F(\lambda)|,
\end{equation*}
and the seminormed space $\mf{L}_{E}^{\infty}(U)\doteqdot\lr{L_{E}^{\infty}(U)}{\|\cdot\|_{\infty}^{E}}$ such that 
\begin{equation*}
\begin{cases}
L_{E}^{\infty}(U)
\coloneqq
\left\{
f:\C\to\C
\mid
f\chi_{U}\in\mr{Bor}(U),\,
|f\chi_{U}|_{\infty}^{E}<\infty
\right\},
\\
(\forall f\in L_{E}^{\infty}(U))
(\|f\|_{\infty}^{E}\coloneqq|f\chi_{U}|_{\infty}^{E}).
\end{cases}
\end{equation*}
Here $\chi_{U}:\C\to\C$ is the characteristic map of $U$ which is by definition equal to $1$ in $U$ and $0$ in $\complement_{\C}U$.
%%%%%%%%%%%%%%%%%%%%%%%%%%%%
\par
Let $Y$ be a locally compact space. In this paper we call Radon measure on $Y$ what Bourbaki in \cite{bourInt}
calls measure on $Y$ whose definition will be here recalled for the sake of completeness. 
More specifically for any $K\in\mr{Comp}(Y)$, with $\mr{Comp}(Y)$ the set of compact subsets
of $Y$, define $\mc{H}(Y,K)$ to be the normed space of complex valued continuous maps defined on $Y$ with compact support 
contained in $K$ and whose norm is $\|f\|_{K}\coloneqq\sup_{y\in K}|f(y)|$, $f\in\mc{H}(Y,K)$. Next define $\mc{H}(Y)$ to be the 
locally convex space of complex valued continuous maps defined on $Y$ with compact support, endowed by the inductive limit of 
the locally convex topologies of the spaces $\{\mc{H}(Y,K)\,\vert\,K\in\mr{Comp}(Y)\}$; 
see \cite[II.29 Proposition 5]{bourTls} for the construction and characterization of final locally convex topologies
and its application in \cite[II.31 exemples II]{bourTls} for the general construction of the inductive limit of locally convex 
spaces. About the topology on $\mc{H}(Y)$ here we need only to know that for any locally convex space $E$, a linear map 
$T:\mc{H}(Y)\to E$ is continuous if and only if $T\circ\imath_{K}$ is continuous for every $K\in\mr{Comp}(Y)$ where 
$\imath_{K}:\mc{H}(Y,K)\to\mc{H}(Y)$ is the map $f\mapsto f$. Equivalently if $\Gamma$ is a fundamental set of seminorms of $E$
then $T$ is continuous if and only if 
\begin{equation*}
(\forall K\in\mr{Comp}(Y))(\forall q\in\Gamma)(\exists\,C>0)(\forall f\in\mc{H}(Y,K))(q(T(f))\leq C\|f\|_{K}).
\end{equation*}
Now a complex valued linear continuous map on $\mc{H}(Y)$ (the above case with $E=\C$) is called measure on $Y$ by 
Bourbaki in \cite[III.7 Definition 2]{bourInt} and \emph{Radon measure on $Y$} by us. Let $\mu$ be a Radon measure on $Y$, 
$\mu$ is called \emph{real} if $\mu(f)\in\R$ for every $f\in\mc{H}(Y;\R)$ with 
$\mc{H}(Y;\R)\coloneqq\{f\in\mc{H}(Y)\,\vert\,(\forall y\in Y)(f(y)\in\R)\}$, while $\mu$ is called \emph{positive} if it is 
real and $\mu(f)\geq 0$ for every $f\in\mc{H}_{+}(Y)$ 
with $\mc{H}_{+}(Y)\coloneqq\{f\in\mc{H}(Y;\R)\,\vert\,(\forall y\in Y)(f(y)\geq 0)\}$.
Let $\mc{S}_{+}(Y)$ be the set of positive lower semi-continuous maps on $Y$ and let $\R_{\geq 0}^{Y}$ be the set of the maps 
defined on $Y$ and with values in $\R_{\geq 0}\coloneqq\{x\in\R\,\vert\,x\geq 0\}$.
If $\nu$ is a positive Radon measure on $Y$, then the \emph{upper integral with respect to $\nu$} is by definition 
the map $\nu^{\ast}:\R_{\geq 0}^{Y}\to\R_{\geq 0}\cup\{+\infty\}$ such that see \cite[IV.6 Definition 3; IV.2 Definition 1]{bourInt}
\begin{equation*}
\begin{cases}
\nu^{\ast}(f)\coloneqq\inf_{h\geq f,h\in\mc{S}_{+}(Y)}\nu^{\ast}(h),\,\forall f\in\R_{\geq 0}^{Y},
\\
\nu^{\ast}(h)\coloneqq\sup_{g\in\mc{H}_{+}(Y),g\leq h}\nu(g),\,\forall h\in\mc{S}_{+}(Y).
\end{cases}
\end{equation*}
Note that the two definitions agree and 
\begin{equation}
\label{11251355}
(\forall f\in\mc{H}_{+}(Y))(\nu^{\ast}(f)=\nu(f)).
\end{equation}
Finally we employ the Bourbaki convention to denote $\mu(f)$ also by $\int_{Y}f(y)d\mu(y)$ with $f\in\mc{H}(Y)$ when $\mu$ is a 
Radon measure on $Y$ and $\nu^{\ast}(g)=\int_{Y}^{\ast}g(y)d\nu(y)$ with $g\in\R_{\geq 0}^{Y}$ when $\nu$ is a positive Radon measure 
on $Y$.
%%%%%%%%%%%%%%%%%%%%%%%%%%%%%%%%%%%%%%%%%%%%%%%%%%%%%  
%%%%%%%%%%%%%%%%%%%%%%%%%%%%%%%%%%%%%%%%%%%%%%%%%%%%%  
\section*{Results}
For the remaining of the paper we fix a complex Banach space $G$, a possibly unbounded scalar type operator $R$ in $G$, 
let $E$ be the resolution of the identity of $R$ and $\sigma(R)$ be the spectrum of $R$.
%%%%%%%%%%%%%%%%%%%%%%%%%%%%%%%%%%%%%%%%%%%%%%%%%%%%%  
%%%%%%%%%%%%%%%%%%%%%%%%%%%%%%%%%%%%%%%%%%%%%%%%%%%%%  
\begin{definition}
Let $F$ denote the locally convex space whose underlying linear space equals $L(\mr{Dom}(R),G)$
and whose locally convex topology is generated by the following set of seminorms 
$\{F\ni P\mapsto \|Pv\|\,\vert\,v\in\mr{Dom}(R)\}$, let $j_{B(G)}^{F}:B(G)\to F$ be the map $T\mapsto T\up\mr{Dom}(R)$.
Let $\widehat{F}$ denote the Hausdorff completion of $F$ and $\mf{i}_{F}$ the canonical map of $F$ into $\widehat{F}$.
\end{definition}
%%%%%%%%%%%%%%%%%%%%%%%%%%%%%%%%%%%%%%%%%%%%%%%%%%%%%  
Notice that $j_{B(G)}^{F}$ and $\mf{i}_{F}$ are linear and continuous maps so uniformly continuous.
%%%%%%%%%%%%%%%%%%%%%%%%%%%%%%%%%%%%%%%%%%%%%%%%%%%%%  
\begin{theorem}
\label{05221753}
Let $K\subseteq\R$ be a compact nondegenerate interval of $\R$, and $U$ be an open neighbourhood of $\sigma(R)$ 
such that 
$K\cdot U\subseteq U$. Furthermore let $X$ be a compact space, 
$\mf{g}\in\mc{C}(X\times K,\mf{L}_{E}^{\infty}(U))$, $\{x_{n}\}_{n\in\N}$ be a sequence in $X$, $t_{0}\in K$ and 
$\{f^{n}\}_{n\in\N}$ be a sequence of analytic maps defined on $U$ with the following properties:
$\sup_{t\in K}\|\widetilde{f_{t}^{n}}\|_{\infty}^{E}<\infty$ for every $n\in\N$,
\begin{equation}
\label{05221753a}
(\forall v\in\mr{Dom}(R))(\exists\,\lim_{n\in\N}f^{n}(t_{0}R)v),
\end{equation}
and for every $n\in\N$ and $t\in K$ we have 
\begin{equation}
\label{05221753b}
\mf{g}(x_{n},t)=\left(\widetilde{\frac{df^{n}}{d\lambda}}\right)_{t}.
\end{equation}
Then there exist $y\in X$ and a subsequence $\{x_{n_{k}}\}_{k\in\N}$ of $\{x_{n}\}_{n\in\N}$ converging to $y$ such that 
for every open neighbourhood $W$ of $t_{0}$ for which $W\subset K$ there exists a map $h\in\mc{C}^{1}(W,\widehat{F})$ whose 
derivative satisfies for every $t\in W$
\begin{equation}
\label{05251055}
(\mr{D}h)(t)=\mf{i}_{F}\bigl(\mf{g}(y,t)(R)R\bigr),
\end{equation}
while
\begin{equation}
\label{05251056}
h=\lim_{k\in\N}(W\ni t\mapsto(\mf{i}_{F}\circ j_{B(G)}^{F})(f^{n_{k}}(tR)))\text{ in }\mc{C}_{u}(W,\widehat{F}).
\end{equation}
\end{theorem}
\begin{proof}
Let $n\in\N$ and let $W$ be an open interval of $\R$ such that $t_{0}\in W$ and $W\subset K$.
We have $\sup_{t\in K}\|\mf{g}(x_{n},t)\|_{\infty}^{E}<\infty$ 
since $\|\cdot\|_{\infty}^{E}\circ\mf{g}$ is continuous and $\{x_{n}\}\times K$ is compact. 
So we can apply \cite[Theorem 1.23]{sil0} and state that for every $v\in\mr{Dom}(R)$ the map 
$W\ni t\mapsto f^{n}(tR)v\in G$ is derivable and for every $t\in W$
\begin{equation}
\label{05222022}
\frac{df^{n}(tR)v}{dt}=R\frac{df^{n}}{d\lambda}(tR)v.
\end{equation}
Now by letting $\imath:\C\to\C$ be the map $\lambda\to\lambda$ clearly we have $R=\imath(R)$ so 
\begin{equation*}
\begin{aligned}
R\frac{df^{n}}{d\lambda}(tR)
&=
R\left(\frac{df^{n}}{d\lambda}\right)_{t}(R)
\\
&=
\imath(R)\left(\frac{df^{n}}{d\lambda}\right)_{t}(R)
\\
&\subseteq
\left(\imath\cdot\left(\frac{df^{n}}{d\lambda}\right)_{t}\right)(R)
\\
&=
\left(\left(\frac{df^{n}}{d\lambda}\right)_{t}\cdot\imath\right)(R);
\end{aligned}
\end{equation*}
where the inclusion follows by the spectral theorem \cite[18.2.11(f)]{ds} which implies also that
\begin{equation*}
\left(\frac{df^{n}}{d\lambda}\right)_{t}(R)\imath(R)
\subseteq
\left(\left(\frac{df^{n}}{d\lambda}\right)_{t}\cdot\imath\right)(R)
\end{equation*}
while clearly 
\begin{equation*}
\mr{Dom}\left(\left(\frac{df^{n}}{d\lambda}\right)_{t}(R)\imath(R)\right)=\mr{Dom}(R);
\end{equation*}
so, for every $v\in\mr{Dom}(R)$ we deduce that 
\begin{equation*}
R\frac{df^{n}}{d\lambda}(tR)v
=
\left(\frac{df^{n}}{d\lambda}\right)_{t}(R)Rv.
\end{equation*}
Therefore by \eqref{05222022} and hypothesis \eqref{05221753b} we obtain for every $v\in\mr{Dom}(R)$ and $t\in W$
\begin{equation}
\label{05222023}
\frac{df^{n}(tR)v}{dt}=\mf{g}(x_{n},t)(R)Rv,
\end{equation}
which is equivalent to state that for every $t\in W$
\begin{equation*}
\frac{dj_{B(G)}^{F}(f^{n}(tR))}{dt}=\mf{g}(x_{n},t)(R)R.
\end{equation*}
%%%%%%%
Next $\mf{g}$ is uniformly continuous being a continuous map defined on a compact set and valued in a uniform space thus, 
by \cite[X.13 Proposition 2]{BourGT} we deduce for every $z\in X$ that 
$Q_{z}\coloneqq\{\mf{g}(z,\cdot)\}\cup\{\mf{g}(x_{n},\cdot)\,\vert\,n\in\N\}$ 
is an equicontinuous subset of $\mc{C}(K,\mf{L}_{E}^{\infty}(U))$ therefore, by \cite[X.16 Theorem 1]{BourGT} on $Q_{z}$ 
the uniform structure inherited by that of $\mc{C}_{u}(K,\mf{L}_{E}^{\infty}(U))$ equals the uniform structure inherited by 
that of $\mc{C}_{s}(K,\mf{L}_{E}^{\infty}(U))$. 
Now let $\{x_{n_{k}}\}_{k\in\N}$ be a subsequence of $\{x_{n}\}_{n\in\N}$ converging to a point $y\in X$ which exists since 
$X$ is compact. Thus, the sequence $\{\mf{g}(x_{n_{k}},\cdot)\}_{k\in\N}$ converges to $\mf{g}(y,\cdot)$ in
$\mc{C}_{s}(K,\mf{L}_{E}^{\infty}(U))$ due to the continuity of $\mf{g}$ therefore, for what above established applied to 
$z=y$ we conclude that
\begin{equation}
\label{06081417}
\mf{g}(y,\cdot)=\lim_{k\in\N}\mf{g}(x_{n_{k}},\cdot)\text{ in }\mc{C}_{u}(K,\mf{L}_{E}^{\infty}(U)).
\end{equation}
Next the functional calculus restricted at $\mf{L}_{E}^{\infty}(U)$ say $\mr{I}$ is a linear and continuous map at values 
in the Banach space $B(G)$ since \cite[18.2.11(c)]{ds} so, uniformly continuous, moreover for every $v\in\mr{Dom}(R)$ 
the map $\hat{v}:B(G)\to G$ defined by $T\mapsto TRv$ is linear and continuous so uniformly continuous, 
therefore, $\hat{v}\circ\mr{I}:\mf{L}_{E}^{\infty}(U)\to G$ is uniformly continuous being composition of uniformly continuous
maps so, we deduce by \eqref{06081417} and by \cite[X.5 Proposition 3(a)]{BourGT} applied to the uniformly continuous map 
$\hat{v}\circ\mr{I}$ that for every $v\in\mr{Dom}(R)$ 
\begin{equation*}
\mf{g}(y,\cdot)(R)Rv=\lim_{k\in\N}\mf{g}(x_{n_{k}},\cdot)(R)Rv\text{ in }\mc{C}_{u}(K,G),
\end{equation*}
and then \eqref{05222023} yields for every $v\in\mr{Dom}(R)$ 
\begin{equation*}
\lim_{k\in\N}
\left(W\ni t\mapsto\frac{df^{n_{k}}(tR)v}{dt}\right)
=
(W\ni t\mapsto\mf{g}(y,t)(R)Rv)
\text{ in }\mc{C}_{u}(W,G),
\end{equation*}
which is equivalent to state that
\begin{equation*}
\lim_{k\in\N}
\left(W\ni t\mapsto\frac{dj_{B(G)}^{F}(f^{n_{k}}(tR))}{dt}\right)
=
(W\ni t\mapsto\mf{g}(y,t)(R)R)
\text{ in }\mc{C}_{u}(W,F),
\end{equation*}
therefore in sequence by \cite[X.5 Proposition 3(a)]{BourGT} and 
the fact that the derivative operator commutes with the operator $f\mapsto u\circ f$ where $u$ is a linear continuous 
operator between Hausdorff locally convex spaces,
both applied to the map $\mf{i}_{F}$ we conclude that 
\begin{equation}
\label{05231256I}
\lim_{k\in\N}\left(W\ni t\mapsto\frac{d(\mf{i}_{F}\circ j_{B(G)}^{F})(f^{n_{k}}(tR))}{dt}\right)
=
(W\ni t\mapsto\mf{i}_{F}(\mf{g}(y,t)(R)R))
\text{ in }\mc{C}_{u}(W,\widehat{F}).
\end{equation}
Next \eqref{05221753a} is equivalent to state that $\exists\,\lim_{n\in\N}j_{B(G)}^{F}(f^{n}(t_{0}R))$ in $F$
which then implies 
\begin{equation}
\label{05231256II}
\exists\,\lim_{n\in\N}(\mf{i}_{F}\circ j_{B(G)}^{F})(f^{n}(t_{0}R))\text{ in }\widehat{F}.
\end{equation}
The statement then follows by \eqref{05231256I}, \eqref{05231256II} and 
the extension of the result \cite[Corollary 8.6.4]{die1} to complete locally convex linear space valued maps defined on open 
sets of $\R$.
\end{proof}
%%%%%%%%%%%%%%%%%%%%%%%%%%%%%%%%%%%%%%%%%%%%%%%%%%  
In what follows we employ the concept of $E$-appropriate set with the duality property which is defined in 
\cite[Definition 2.11]{sil0}.
%%%%%%%%%%%%%%%%%%%%%%%%%%%%%%%%%%%%%%%%%%%%%%%%%%  
\begin{definition}
Let $Y$ be a locally compact space, $\mu$ be a positive Radon measure on $Y$, $X$ be a set, $U$ be an open neighbourhood of 
$\sigma(R)$, $\mf{h}:X\times Y\to\mf{L}_{E}^{\infty}(U)$, $\{z_{k}\}_{k\in\N}$ be a sequence in $X$, and $\mc{N}$ be 
an $E$-appropriate set with the duality property. We say that $\lr{\mf{h}}{\{z_{k}\}_{k\in\N}}$ satisfies the 
\emph{bounding assumption with respect to $\mc{N}$} iff for every $k\in\N$ 
\begin{itemize}
\item
\begin{equation}
\label{11242149}
\int_{Y}^{\ast}\|\mf{h}(z_{k},y)\|_{\infty}^{E}d\mu(y)<\infty;
\end{equation}
\item
for all $\omega\in\mc{N}$ the map $Y\ni y\mapsto\omega(\mf{h}(z_{k},y)(R))\in\C$ is $\mu$-measurable.
\end{itemize}
\end{definition}
%%%%%%%%%%%%%%%%%%%%%%%%%%%%%%%%%%%%%%%%%%%%%%%%%%  
\begin{lemma}
\label{11242141}
Let $Y$ be a locally compact space, $\mu$ be a positive Radon measure on $Y$, $X$ be a topological space, $U$ be an open 
neighbourhood of $\sigma(R)$, $\{z_{k}\}_{k\in\N}$ be a sequence in $X$, and $\mc{N}$ be an $E$-appropriate set with 
the duality property. If $\mf{h}\in\mc{C}(X\times Y,\mf{L}_{E}^{\infty}(U))$ and for every $k\in\N$ the map 
$Y\ni y\mapsto\mf{h}(z_{k},y)\in\mf{L}_{E}^{\infty}(U)$ has compact support, then $\lr{\mf{h}}{\{z_{k}\}_{k\in\N}}$ satisfies the 
bounding assumption with respect to $\mc{N}$.
\end{lemma}
\begin{proof}
Let $k\in\N$. The map $Y\ni y\mapsto\|\mf{h}(z_{k},y)\|_{\infty}^{E}$ has compact support by hypothesis and it is continuous 
since are continuous the map $Y\ni y\mapsto\mf{h}(z_{k},y)\in\mf{L}_{E}^{\infty}(U)$ and the seminorm $\|\cdot\|_{\infty}^{E}$, 
then \eqref{11242149} follows since \eqref{11251355}. Next the restriction $\mf{L}_{E}^{\infty}(U)\ni l\mapsto l(R)\in B(G)$ of the
functional calculus of $R$ at $\mf{L}_{E}^{\infty}(U)$ is continuous being bounded since \cite[Theorem 18.2.11(c)]{ds} then for 
all $\omega\in\mc{N}$ the map $Y\ni y\mapsto\omega(\mf{h}(z_{k},y)(R))\in\C$ is $\mu$-measurable being continuous.
\end{proof}
%%%%%%%%%%%%%%%%%%%%%%%%%%%%%%%%%%%%%%%%%%%%%%%%%%  
Thm.\ref{05221753} combined with our \cite[Corollary 2.33]{sil0} yields the following 
%%%%%%%%%%%%%%%%%%%%%%%%%%%%%%%%%%%%%%%%%%%%%%%%%%  
\begin{corollary}
\label{11251048}
Let the hypothesis of Thm.\ref{05221753} be satisfied so, there exists $y\in X$ and $\{x_{n_{k}}\}_{k\in\N}$ a
subsequence of $\{x_{n}\}_{n\in\N}$ converging to $y$.
Furthermore let $\mc{N}$ be an $E$-appropriate set with the duality property. 
Then for every open interval $W$ such that $t_{0}\in W$ and 
$W\subseteq\overset{\circ}{K}$, and for every $u_{1},u_{2}\in W$ we have 
\begin{equation*}
\int_{u_{1}}^{u_{2}}\mf{i}_{F}\bigl(\mf{g}(y,t)(R)R\bigr)dt
=
\lim_{k\in\N}(\mf{i}_{F}\circ j_{B(G)}^{F})\bigl(R\int_{u_{1}}^{u_{2}}\mf{g}(x_{n_{k}},t)(R)dt\bigr).
\end{equation*}
Here the limit is in the space $\widehat{F}$, the integral in the left-hand side is the weak integral 
of the map $[u_{1},u_{2}]\mapsto t\ni\mf{i}_{F}\bigl(\mf{g}(y,t)(R)R\bigr)\in\widehat{F}$
with respect to the Lebesgue measure on $[u_{1},u_{2}]$ and the $\sigma(\widehat{F},\widehat{F}^{\ast})$-topology
with $\widehat{F}^{\ast}$ the topological dual of $\widehat{F}$; while the integral in the right-hand side is the weak integral 
of the map $[u_{1},u_{2}]\ni t\mapsto\mf{g}(x_{n_{k}},t)(R)\in B(G)$ with respect to the Lebesgue measure on $[u_{1},u_{2}]$ and 
the $\sigma(B(G),\mc{N})$-topology.
\end{corollary}
\begin{proof}
Let $h$ be the map of which in the statement of Thm.\ref{05221753} thus, 
\begin{equation*}
\begin{aligned}
\int_{u_{1}}^{u_{2}}\mf{i}_{F}\bigl(\mf{g}(y,t)(R)R\bigr)dt
&=
h(u_{2})-h(u_{1})
\\
&=
\lim_{k\in\N}(\mf{i}_{F}\circ j_{B(G)}^{F})
(f^{n_{k}}(u_{2}R)-f^{n_{k}}(u_{1}R))
\\
&=
\lim_{k\in\N}(\mf{i}_{F}\circ j_{B(G)}^{F})
\bigl(R\int_{u_{1}}^{u_{2}}\frac{df^{n_{k}}}{d\lambda}(tR)dt\bigr)
\\
&=
\lim_{k\in\N}(\mf{i}_{F}\circ j_{B(G)}^{F})
\bigl(R\int_{u_{1}}^{u_{2}}\mf{g}(x_{n_{k}},t)(R)dt\bigr);
\end{aligned}
\end{equation*}
where the first equality follows by equality \eqref{05251055} in Theorem \ref{05221753} 
and by the fundamental theorem of calculus for locally convex space valued maps (see for instance \cite[Theorem 1.5]{glo});
the second equality follows from the construction of the map $h$ in \eqref{05251056};
the third equality follows since \cite[Corollary 2.33]{sil0} and since $\lr{\mf{g}}{\{x_{n_{k}}\}_{k\in\N}}$ satisfies the bounding 
assumption with respect to $\mc{N}$ by Lemma \ref{11242141}; the fourth equality follows since hypothesis \eqref{05221753b}.
\end{proof}
%%%%%%%%%%%%%%%%%%%%%%%%%%%%%%%%%%%%%%%%%%%%%%%%%%  
\section*{Data availability statement and declaration}
The author declares that the data supporting the findings of this study are available within the paper.
Furthermore the author has no competing interests to declare that are relevant to the content of this article.
%%%%%%%%%%%%%%%%%%%%%%%%%%%%%%%%%%%%%%%%%%%%%%%%%%  


\begin{thebibliography}{99}

\bibitem{BourGT}
N. Bourbaki, 
\emph{Topologie Generale}, CCLS, Ch 5-10 1974.


\bibitem{bourTls}
N. Bourbaki, 
\emph{Espaces Vectoriels Topologiques}
Masson 1981

\bibitem{bourInt}
N.Bourbaki,
\emph{Integration I. Chapters 1-6}
Springer 2004


\bibitem{die1}
J. Dieudonne,
\emph{Foundations of Modern Analysis}
Academic Press 1969.

\bibitem{ds}
N. Dunford; J. T. Schwartz,
\emph{Linear Operators Part III}
Wiley Interscience Publ. 1971.

\bibitem{glo}
H. Gl\"{o}ckner,
\emph{Infinite-dimensional Lie groups without completeness restrictions} 
In: Geometry and analysis on finite- and infinite-dimensional Lie groups (B\c{e}dlewo,2000)
Banach Center Publ., 55, Polish Acad. Sci. Inst. Math., Warsaw, 2002, pp. 43–59.


\bibitem{sil0}
B. Silvestri,
\emph{Integral equalities for functions of unbounded spectral operators in Banach spaces}
Dissertationes Math. 464 (2009), 60 pp.

\end{thebibliography}
\end{document}